\documentclass[number,citesort,seceqn,dvips]{arxbj}

\aid{0}
\volume{17}
\issue{3}
\pubyear{2011}
\firstpage{1044}
\lastpage{1053}
\doi{10.3150/10-BEJ302}

\makeatletter

\newtheorem{theorem}{Theorem}
\newtheorem{corollary}{Corollary}
\newtheorem{lemma}{Lemma}
\newtheorem{proposition}{Proposition}
\newtheorem{claim}{Claim}

\newremark{remark}{Remark}

\makeatother

\begin{document}
\begin{frontmatter}

\title{Some stochastic inequalities for weighted sums}
\runtitle{Stochastic inequalities for weighted sums}

\begin{aug}
\author{\fnms{Yaming} \snm{Yu}\corref{}\ead[label=e1]{yamingy@uci.edu}}
\runauthor{Y. Yu}
\address{Department of Statistics,
University of California,
Irvine, CA 92697-1250, USA.\\
\printead{e1}}
\end{aug}

\received{\smonth{10} \syear{2009}}

%
\begin{abstract}
We compare weighted sums of i.i.d. positive random variables according
to the usual stochastic order. The
main inequalities are derived using majorization techniques under
certain log-concavity
assumptions. Specifically, let $Y_i$ be i.i.d. random variables on
$\mathbf{R}_+$. Assuming that $\log Y_i$
has a log-concave density, we show that $\sum a_i Y_i$ is
stochastically smaller than $\sum b_i Y_i$,
if $(\log a_1, \ldots, \log a_n)$ is majorized by $(\log b_1, \ldots,
\log b_n)$. On the other hand, assuming that
$Y_i^p$ has a log-concave density for some $p>1$, we show that $\sum
a_i Y_i$ is stochastically larger than
$\sum b_i Y_i$, if $(a_1^q,\ldots, a_n^q)$ is majorized by $(b_1^q,
\ldots, b_n^q)$, where $p^{-1}+q^{-1}=1$.
These unify several stochastic ordering results for specific
distributions. In particular, a~conjecture of Hitczenko [\textit
{Sankhy\={a} A} \textbf{60}
(1998) 171--175] on Weibull variables is proved. Potential applications
in reliability and wireless communications are mentioned.
\end{abstract}

%

\end{frontmatter}

\section{Main results and examples}

This paper aims to unify and generalize certain stochastic comparison
results concerning weighted sums. Let
$Y_1,\ldots, Y_n$ be i.i.d. random variables on $\mathbf{R}_+$. We are
interested in comparing two weighted sums, $\sum_{i=1}^n a_i Y_i$ and
$\sum_{i=1}^n b_i Y_i, a_i, b_i\in\mathbf{R}_+,$ with respect to the
usual stochastic order. A random variable $X$ is said to be no larger
than $Y$ in the usual stochastic order, written as $X\leq_{\mathrm{st}} Y$, if
$\Pr(X>t)\leq\Pr(Y>t)$ for all $t\in\mathbf{R}$. For an introduction
to various stochastic orders, see \cite{2007Shaked}.
Ordering in terms of $\leq_{\mathrm{st}}$ may be used to bound the tail
probability of $\sum a_iY_i,$ for example, in terms of the tail
probability of $\sum Y_i$. For specific distributions, such comparisons
have been explored in several contexts, including reliability \cite{1994Boland,1999Bon}.

We shall use the notion of \textit{majorization} \cite{1979Marshall}. A real vector $b=(b_1, \ldots, b_n)$ is said to majorize
$a=(a_1, \ldots, a_n)$, written as $a\prec b$, if (i) $\sum_{i=1}^n
a_i=\sum_{i=1}^n b_i$, and (ii) $\sum_{i=k}^n a_{(i)}\leq\sum_{i=k}^n
b_{(i)}, k=2, \ldots, n,$ where $a_{(1)}\leq\cdots\leq a_{(n)}$ and
$b_{(1)}\leq\cdots\leq b_{(n)}$ are $(a_1, \ldots, a_n)$ and $(b_1,
\ldots, b_n)$ arranged in increasing order, respectively. A function
$\phi(a)$ symmetric in the coordinates of $a=(a_1, \ldots, a_n)$ is
said to be \textit{Schur-concave} if $a\prec b$ implies $\phi(a)\geq
\phi(b).$ A~function $\phi(a)$ is \textit{Schur-convex} if $-\phi(a)$
is Schur-concave.

A non-negative function $f(x), x\in\mathbf{R}^n,$ is \textit
{log-concave} if $\operatorname{supp}(f)$ is convex and $\log f(x)$ is concave on
$\operatorname{supp}(f)$. Log-concavity plays a critical role in deriving our main
results. For other stochastic comparison results involving
log-concavity, see, for example, \cite{1981Karlin,2008Yu,2009bYu,2010Yu}.

In this section, after stating our main results (Theorems \ref{thm1}
and \ref{thm2}), we illustrate with several examples and mention
potential applications. The main results are proved in Section \ref{sec2}. Some
technical details in the proof of Theorem \ref{thm2} are collected in
the \hyperref[app]{Appendix}.

\begin{theorem}
\label{thm1}
Let $Y_1, \ldots, Y_n$ be i.i.d. random variables with density $f(y)$
on $\mathbf{R}_+$ such that $f(\mathrm{e}^x)$ is log-concave in $x\in\mathbf
{R}$. Then, for each $t>0, \Pr(\sum a_i Y_i\leq t)$ is a Schur-concave
function of
\[
\log a\equiv(\log a_1, \ldots, \log a_n).
\]
Equivalently, if $a, b\in\mathbf{R}_+^n$, then
%
\begin{equation}
\label{ineqthm1}
\log a\prec\log b\quad \Longrightarrow\quad \sum a_i Y_i\leq_{\mathrm{st}} \sum b_i Y_i.
\end{equation}
\end{theorem}

\begin{theorem}
\label{thm2}
Let $p>1,$ and let $Y_1,\ldots, Y_n$ be i.i.d. random variables with
density $f(y)$ on $\mathbf{R}_+$ such that the function
%
\begin{equation}
\label{lc}
\min\{0, 2/p-1\}\log x + \log f(x^{1/p})
\end{equation}
is concave in $x\in\mathbf{R}_+$. Then, for each $t>0, \Pr(\sum a_i
Y_i \leq t)$ is a Schur-convex function of
\[
a^q\equiv(a_1^q, \ldots, a_n^q)\in\mathbf{R}_+^n,
\]
where $p^{-1}+q^{-1}=1$. Equivalently, if $a, b\in\mathbf{R}_+^n$, then
%
\begin{equation}
\label{ineqthm2}
a^q\prec b^q \quad\Longrightarrow\quad \sum b_i Y_i \leq_{\mathrm{st}} \sum a_i Y_i.
\end{equation}
\end{theorem}

\begin{remark*} In Theorem \ref{thm1}, the condition that $f(\mathrm{e}^x)$ is
log-concave is equivalent to $\log Y_i$ having a log-concave density
(see, e.g., \cite{2009Righter}). In Theorem \ref{thm2}, a
sufficient condition for (\ref{lc}) is that $x^{1/p-1}f(x^{1/p})$ is
log-concave, or, equivalently, $Y_i^{p}$ has a log-concave density
(this special case is mentioned in the abstract). Theorems \ref{thm1}
and \ref{thm2} are quite applicable, since log-concavity is associated
with many well-known densities (see Corollaries \ref{coro1} and \ref
{coro2}).

Theorem \ref{thm1} is reminiscent of the following result of \cite{1965Proschan},
originally stated in terms of the peakedness order.
\end{remark*}

\begin{theorem}
\label{thm4}
Let $Y_i, i=1,\ldots,n,$ be i.i.d. random variables on $\mathbf{R}$
with a log-concave density that is symmetric about zero. Then for each
$t>0, \Pr(\sum a_i Y_i\leq t)$ is a Schur-concave function of $a\in
\mathbf{R}_+^n$.
\end{theorem}

Theorem \ref{thm2} is closely related to Theorem \ref{kr}, which is a
version (with a stronger assumption) of Theorem 24 of \cite{1983Karlin}.

\begin{theorem}
\label{kr}
Let $0<p<1,$ and let $Y_1,\ldots, Y_n$ be i.i.d. random variables on
$\mathbf{R}_+$ such that $Y_i^p$ has a log-concave density. Then, for
each $t>0, \Pr(\sum a_i Y_i \leq t)$ is a Schur-concave function of
$(a_1^q, \ldots, a_n^q)\in\mathbf{R}_+^n,$ where $p^{-1}+q^{-1}=1$.
\end{theorem}

Karlin and Rinott \cite{1983Karlin} gave an elegant proof of Theorem \ref{kr}
using the Pr\'{e}kopa--Leindler inequality. Our proofs of Theorems \ref
{thm1} and \ref{thm2} (Section \ref{sec2}) borrow ideas from both \cite{1965Proschan}
 and~\cite{1983Karlin}. See \cite{1988Shaked} for
more related inequalities.

Bounds on the distribution function of $\sum a_i Y_i$ are readily
obtained in terms of the distribution function of $\sum Y_i$. In
Theorem \ref{thm1}, for example, (\ref{ineqthm1}) gives
%
\begin{equation}
\label{upper}
\Pr\Bigl(\sum b_i Y_i \leq t\Bigr)\leq\Pr\Bigl(b_*\sum Y_i\leq t\Bigr),\qquad b_*=\Bigl(\prod
b_i\Bigr)^{1/n}, t>0.
\end{equation}
In Theorem \ref{thm2}, (\ref{ineqthm2}) gives
%
\begin{equation}
\label{lower}
\Pr\Bigl(\sum b_i Y_i \leq t\Bigr)\geq\Pr\Bigl(b^* \sum Y_i\leq t\Bigr),\qquad
b^*=\Bigl(n^{-1}\sum b_i^q\Bigr)^{1/q}, t>0.
\end{equation}
More generally, we obtain inequalities for the expectations of monotone
functions, since $X\leq_{\mathrm{st}} Y$ implies
${Eg}(X)\leq {Eg}(Y)$ for every increasing function $g$ such that the
expectations exist.

Let us mention some specific distributions to which Theorems \ref{thm1}
and \ref{thm2} can be applied. Corollary \ref{coro1} follows from
Theorem \ref{thm1}. The log-concavity condition is easily verified in
each case (for more distributions that satisfy this condition, see \cite{2004Hu},
Example 1). Related results on sums of uniform
variables can be found in \cite{2002Korwar}. The gamma case has recently
been discussed by Khaledi and Kochar \cite{2004Khaledi}, Yu \cite{2009aYu} and Zhao and
Balakrishnan \cite{2009Zhao}.

\begin{corollary}
\label{coro1}
For $a, b\in\mathbf{R}_+^n$, (\ref{ineqthm1}) holds when $Y_i$ are
i.i.d. having one of the following distributions:
\begin{enumerate}[(5)]
\item[(1)]
uniform on the interval $(0, s), s>0$;
\item[(2)]
gamma$(\alpha, \beta), \alpha, \beta>0$;
\item[(3)]
any log-normal distribution;
\item[(4)]
the Weibull distribution with parameter $p>0$, whose density is
\[
f(y)=p y^{p-1} \mathrm{e}^{-y^p},\qquad y>0;
\]
\item[(5)]
the generalized Rayleigh distribution with parameter $\nu>0,$ whose
density is
\[
f(y)\propto y^{\nu-1} \mathrm{e}^{-y^2/2},\qquad y>0.
\]
\end{enumerate}
\end{corollary}

The inequality (\ref{upper}) holds for each of these distributions. The
gamma case is interesting in that the upper bound in (\ref{upper}) is
in terms of a single gamma variable, $\sum Y_i$. The gamma case with
$\alpha=1/2$ dates back to \cite{1960Okamoto}. See also \cite{1987Bock,2003Sz} for related inequalities.

\begin{corollary}
\label{coro2}
Let $p>1,$ and define $q$ by $p^{-1}+q^{-1}=1$. Then, for $a, b\in
\mathbf{R}_+^n$, (\ref{ineqthm2})
holds in the following cases:
\begin{enumerate}[(1)]
\item[(1)]
$Y_i$ are i.i.d. Weibull variables with parameter $p$;
\item[(2)]
$p=q=2$, and $Y_i$ are i.i.d. generalized Rayleigh variables with
parameter $\nu\geq1.$
\end{enumerate}
\end{corollary}

Corollary \ref{coro2} follows from Theorem \ref{thm2}. The condition
(\ref{lc}) is easily verified. For example, in case 1, $Y_i^p$ has a
log-concave density, which implies (\ref{lc}). Case 1 confirms a
conjecture of \cite{1998Hitczenko}. Case 2 recovers some results of \cite{2000Hu,2001Hu}.

The Weibull case and the generalized Rayleigh case are interesting in
that Corollary \ref{coro1} is also applicable, and we obtain a double
bound through (\ref{upper}) and (\ref{lower}). For example, if $Y_i$
are i.i.d. generalized Rayleigh variables with parameter $\nu\geq1$, then
\[
\Pr\Bigl(a^*\sum Y_i\leq t\Bigr)\leq\Pr\Bigl(\sum a_i Y_i \leq t\Bigr)\leq\Pr\Bigl( a_* \sum
Y_i\leq t\Bigr),\qquad a_i>0, t>0,
\]
where $a^*=(n^{-1}\sum a_i^2)^{1/2}$ and $a_*=(\prod a_i)^{1/n}.$
Manesh and Khaledi \cite{2008Manesh} present related inequalities.

We briefly mention some applications:
\begin{itemize}
\item
Weighted sums of independent $\chi^2$ variables arise naturally in
multivariate statistics as quadratic forms in normal variables.
Stochastic comparisons between such weighted sums are therefore
statistically interesting, and can lead to bounds on the distribution functions.
\item
Suppose the component lifetimes of a redundant standby system (without
repairing) are modeled by a scale family of distributions. Then the
total lifetime is of the form $\sum_i a_i Y_i$. When $Y_i$ are i.i.d.
exponential variables, Bon and Paltanea \cite{1999Bon} obtain comparisons of
the total lifetime with respect to several stochastic orders. Our
Corollary \ref{coro1} shows that, for the~usual stochastic order, (\ref
{ineqthm1}) actually holds for a broad class of distributions,
including the commonly used gamma, Weibull, and log-normal distributions.
\item
When $Y_i$ are i.i.d. exponential variables and $a_i\in\mathbf{R}_+$,
the quantity $E\log(1+\sum a_i Y_i)$ appears in certain wireless
communications problems \cite{2007Jorswieck}. By the
monotonicity of $\log(1+x)$, we have
\[
\sum a_i Y_i\leq_{\mathrm{st}} \sum b_i Y_i\quad \Longrightarrow\quad E\log\Bigl(1+\sum a_i
Y_i\Bigr)\leq E\log\Bigl(1+\sum b_i Y_i\Bigr).
\]
Corollary \ref{coro1} therefore leads to qualitative comparisons for
this expected value. Other weighted sums (e.g., of Rayleigh variables)
also appear in the context of communications.
\end{itemize}

It would be interesting to see whether results similar to Theorems \ref
{thm1}, \ref{thm2} and \ref{kr}
can be obtained for the hazard rate order, or the likelihood ratio
order. For sums of independent gamma variables, such results have been
obtained by Boland, El-Neweihi and~Pro\-schan~\cite{1994Boland}, Bon and Paltanea \cite{1999Bon}, Korwar
\cite{2002Korwar}, Khaledi and Kochar (2004) \cite{2004Khaledi}, Yu (2009)~\cite{2009aYu} and Zhao and
Balakrishnan (2009) \cite{2009Zhao}.

\section{Proofs}\label{sec2}

Two proofs are presented for Theorem \ref{thm1}. The first one uses the
Pr\'{e}kopa--Leindler inequality (Lemma \ref{lem1}) and is inspired by
Karlin and Rinott \cite{1983Karlin}.

\begin{lemma}
\label{lem1}
If $g(x, y)$ is log-concave in $(x, y)\in\mathbf{R}^m\times\mathbf
{R}^n,$ then $\int_{\mathbf{R}^m} g(x, y)\, {\rm d}x$ is log-concave in
$y\in\mathbf{R}^n$.
\end{lemma}

We also use a basic criterion for Schur-concavity.
\begin{proposition}
\label{schur}
If $h(\alpha), \alpha\in\mathbf{R}^n,$ is log-concave and permutation
invariant in $\alpha$, then it is Schur-concave.
\end{proposition}

\begin{pf*}{First proof of Theorem \ref{thm1}}
For $t>0$, define
\[
g(x, \alpha)\equiv1_K \prod_{i=1}^n \mathrm{e}^{x_i}f(\mathrm{e}^{x_i}),\qquad K\equiv\Biggl\{
(x, \alpha)\in\mathbf{R}^{2n}\dvtx \sum_{i=1}^n \mathrm{e}^{x_i+\alpha_i}\leq t\Biggr\},
\]
where $x=(x_1, \ldots, x_n), \alpha=(\alpha_1, \ldots, \alpha_n)\in
\mathbf{R}^n$. Note that $K$ is a convex set ($1_K$ denotes the
indicator function). Since $f(\mathrm{e}^{x_i})$ is log-concave, we know that
$g(x, \alpha)$ is log-concave in $(x, \alpha)$. By Lemma \ref{lem1},
\[
h(\alpha)\equiv\Pr\Bigl(\sum \mathrm{e}^{\alpha_i}Y_i\leq t\Bigr)=\int_{\mathbf{R}^n}
g(x, \alpha)\, \mathrm{d}x
\]
is log-concave in $\alpha\in\mathbf{R}^n$. Since $h(\alpha)$ is
permutation invariant, it is Schur-concave in $\alpha$ by Proposition
\ref{schur}, and the claim is proved.
\end{pf*}

The second proof is inspired by Proschan \cite{1965Proschan}, and serves as an
introduction to the proof of Theorem \ref{thm2}.
Properties of majorization imply that it suffices to prove (\ref
{ineqthm1}) for $a\prec b$ such that $a$ and $b$ differ only in two
components. Since $\leq_{\mathrm{st}}$ is closed under~convolu\-tion~\cite{2007Shaked}, we only need to prove (\ref{ineqthm1}) for $n=2$.

We shall use log-concavity in the following form. If $g(x), x\in
\mathbf{R},$ is log-concave, and $(x_1, x_2)\prec(y_1, y_2)$, then
\[
g(x_1)g(x_2)-g(y_1)g(y_2)\geq0.
\]

\begin{pf*}{Second proof of Theorem \ref{thm1}}
Fix $t>0$, and let $F$ denote the distribution function of~$Y_1$. It
suffices to show that
\[
h(\beta)\equiv\Pr(\beta^{-1} Y_1+ \beta Y_2 \leq t)=\int_0^\infty
F(t\beta-\beta^2 y) f(y) \,\mathrm{d} y
\]
increases in $\beta\in(0, 1]$. We may assume that $\operatorname{supp}(f)\subset
[\varepsilon, \infty)$ for some $\varepsilon>0$. The general case
follows by a standard limiting argument. We can then justify
differentiation under the integral sign and obtain
\begin{eqnarray}\label{inte}
h'(\beta) &=&\int_0^\infty(t-2\beta y) f(t\beta- \beta^2 y) f(y) \,\mathrm{d}y\nonumber
\\[-8pt]\\[-8pt]
&=&\int_0^{t/(2\beta)} (t-2\beta y) f(t\beta-\beta^2 y) f(y) \,\mathrm{d}y
+\int_{t/(2\beta)}^{t/\beta} (t-2\beta y) f(t\beta- \beta^2 y) f(y)\nonumber
\,\mathrm{d}y.\quad
\end{eqnarray}
%

By a change of variables $y\to t/\beta-y$ in the second integral in
(\ref{inte}), we get
\[
h'(\beta)=\int_0^{t/(2\beta)} (t-2\beta y) [f(t\beta-\beta^2 y)f(y)
-f(\beta^2 y) f(t/\beta-y)] \,\mathrm{d}y.
\]
If $0<y<t/(2\beta)$ and $0<\beta\leq1$, then $\beta^2 y\leq\min\{y,
t\beta-\beta^2 y\}.$ That is,
\[
\bigl(\log(t\beta-\beta^2 y), \log y\bigr) \prec\bigl(\log(\beta^2 y), \log(t/\beta-y)\bigr).
\]
Since $f(\mathrm{e}^x)$ is log-concave, we have
\[
f(t\beta-\beta^2 y)f(y) -f(\beta^2 y) f(t/\beta-y)\geq0,\qquad
0<y<t/(2\beta),
\]
which leads to $h'(\beta)\geq0$, as required.
\end{pf*}

Our proof of Theorem \ref{thm2} is similar to (but more involved than)
the second proof of Theorem \ref{thm1}. Under the stronger assumption
that $Y_i^p$ has a log-concave density, we actually obtain a simpler
proof of Theorem \ref{thm2} following the first proof of Theorem \ref
{thm1} (see \cite{1983Karlin}). It seems difficult, however, to
extend this argument assuming only that (\ref{lc}) is concave.

\begin{pf*}{Proof of Theorem \ref{thm2}}
We may assume $n=2$ as in the second proof of Theorem \ref{thm1}. Fix
$t>0$. Effectively, we need to show that
\[
h(\beta)\equiv\Pr\bigl(\beta^{1/q} Y_1+ (1-\beta)^{1/q} Y_2 \leq t\bigr)=\int
_0^\infty F\bigl(t\beta^{-1/q}-(\beta^{-1}-1)^{1/q} y\bigr) f(y) \,\mathrm{d} y
\]
increases in $\beta\in[1/2, 1)$ ($F$ denotes the distribution function
of $Y_1$). We have
\begin{eqnarray*}
q(1-\beta)^{1/p}\beta^{1/q+1} h'(\beta) =\int_0^\infty g(y) \,\mathrm{d}y
=\int_0^{y_0} g(y) \,\mathrm{d}y +\int_{y_0}^{y_1} g(y) \,\mathrm{d}y,
\end{eqnarray*}
where
%
\begin{equation}
\label{lims}
y_0=t(1-\beta)^{1/p},\qquad y_1=t(1-\beta)^{-1/q}
\end{equation}
and
%
\begin{equation}
\label{gxy}
g(y) =(y-y_0) f(x(y))f(y),\qquad x(y) = (\beta^{-1}-1)^{1/q}(y_1-y).
\end{equation}
Differentiation under the integral sign is permitted because
\[
|g(y)|\leq(y_1-y_0) M f(y),\qquad0<y<y_1,
\]
where $M=\sup_{y>0} f(y)$. We know $M<\infty$ because (\ref{lc})
implies that $f(x^{1/p})$ is log-concave in $x\in\mathbf{R}_+$.

In the \hyperref[app]{Appendix}, we prove:

\begin{claim}
\label{claim1}
For each $y\in(0, y_0)$, there exists a unique $\tilde{y}\in(y_0,
y_1)$ such that
%
\begin{equation}
\label{map}
y^p+x^p(y)=\tilde{y}^p+x^p(\tilde{y}),
\end{equation}
where $x(y)$ is given by (\ref{gxy}).
\end{claim}

Henceforth let $y$ and $\tilde{y}$ be related by (\ref{map}). Direct
calculation using the implicit function theorem gives
\begin{eqnarray*}
\frac{\mathrm{d}\tilde{y}}{\mathrm{d} y} =\frac{(\beta y)^{p/q}-((1-\beta
)(y_1-y))^{p/q}}{(\beta\tilde{y})^{p/q}-((1-\beta)(y_1-\tilde{y}))^{p/q}}.
\end{eqnarray*}
A change of variables $y\to\tilde{y}$ in $\int_0^{y_0} g(y) \,\mathrm{d}y$ yields
\[
q(1-\beta)^{1/p}\beta^{1/q+1} h'(\beta)=\int_A g(y) \biggl|\frac{\mathrm{d}y}{\mathrm{d}\tilde{y}}\biggr| \,\mathrm{d}\tilde{y}
+ \int_{y_0}^{y_1} g(z) \,\mathrm{d}z,
\]
where $A\subset(y_0, y_1)$ is the image of the interval $(0, y_0)$
under the mapping $y\to\tilde{y}$. Note that $g(z)\geq0$ for $y_0< z
<y_1$. Hence
\begin{eqnarray}\label{hprime}
q(1-\beta)^{1/p}\beta^{1/q+1} h'(\beta) &\geq&\int_A \biggl(g(\tilde{y})+
g(y) \biggl|\frac{\mathrm{d}y}{\mathrm{d}\tilde{y}}\biggr|\biggr) \,\mathrm{d}\tilde{y}\nonumber
\\[-8pt]\\[-8pt]
&\geq&\int_A (\tilde{y}-y_0) \biggl[f(x(\tilde{y}))f(\tilde{y})-\biggl( \frac{x(y)
y}{x(\tilde{y})\tilde{y}}\biggr)^\delta f(x(y))f(y)\biggr] \,\mathrm{d}\tilde{y},\qquad\ \nonumber
\end{eqnarray}
where $\delta=\min\{0, 2-p\}$. The inequality (\ref{hprime}) is deduced
from Claim \ref{claim2}, which we prove in the \hyperref[app]{Appendix}.

\begin{claim}
\label{claim2}
We have
%
\begin{equation}
\label{ineq3}
\biggl|\frac{\mathrm{d}\tilde{y}}{\mathrm{d}y}\biggr|\geq\biggl( \frac{x(\tilde{y}) \tilde
{y}}{x(y)y}\biggr)^\delta\biggl(\frac{y_0 -y}{\tilde{y}-y_0}\biggr),\qquad0<y<y_0.
\end{equation}
\end{claim}

In the \hyperref[app]{Appendix} we also show:

\begin{claim}
\label{claim3}
For $0<y<y_0,$ we have\vspace*{-2pt}
%
\begin{eqnarray}
\label{uv1}
\beta\tilde{y} &\geq&(1-\beta)(y_1-y);
\\
\label{uv2}
\beta y &\leq&(1-\beta)(y_1-\tilde{y}).
\end{eqnarray}
\end{claim}

For $0<y<y_0$, (\ref{uv2}) yields $y\leq\min\{\tilde{y}, y_1-\tilde
{y}\}$, that is,
\[
(\tilde{y}, y_1-\tilde{y})\prec(y, y_1-y).
\]
Thus, $y(y_1-y)\leq\tilde{y}(y_1-\tilde{y}),$ or, equivalently, $y^p
x^p(y)\leq\tilde{y}^p x^p(\tilde{y}).$ By (\ref{map}), this implies
the relation\vspace*{-2pt}
\[
(\tilde{y}^p, x^p(\tilde{y}))\prec(y^p, x^p(y)),\qquad0<y<y_0.
\]
Assumption (\ref{lc}) then yields ($\delta=\min\{0, 2-p\}$)
\[
(x(\tilde{y})\tilde{y})^\delta f(x(\tilde{y}))f(\tilde{y})-(x(y)
y)^\delta f(x(y))f(y)\geq0,\qquad\tilde{y}\in A.
\]
It follows that the integrand in (\ref{hprime}) is non-negative, and
$h'(\beta)\geq0, \beta\in[1/2, 1),$ as required.
\end{pf*}

\begin{remark} The main complication in the proof of Theorem \ref{thm2}
is that the mapping $y\to\tilde{y}$ is not in closed form. In the
special case $p=q=2,$ where $\tilde{y}$ is explicitly available, the
proof can be simpler.
\end{remark}

\begin{appendix}

\section*{\texorpdfstring{Appendix: Proofs of Claims \protect\ref{claim1}--\protect\ref{claim3}}%
{Appendix: Proofs of Claims 1--3}}\label{app}

It is convenient to prove Claims \ref{claim1}, \ref{claim3} and \ref{claim2} in that order. We emphasize
that no circular argument is involved.
\begin{pf*}{Proof of Claim \ref{claim1}}
Define
\setcounter{equation}{0}
\begin{equation}
\label{ly}
L(y)=\beta^{p/q} y^p +(1-\beta)^{p/q}(y_1-y)^p,\qquad 0\leq y\leq y_1,
\end{equation}
where $y_1$ is given by (\ref{lims}). We have
\[
L'(y)=p\beta^{p/q} y^{p-1}-p(1-\beta)^{p/q} (y_1-y)^{p-1},
\]
and the unique solution of $L'(y)=0$ is $y_0=t(1-\beta)^{1/p}$. Moreover,
\[
L''(y)=p(p-1)[\beta^{p/q} y^{p-2}+(1-\beta)^{p/q} (y_1-y)^{p-2}]> 0.
\]
Hence $L(y)$ strictly decreases on the interval $(0, y_0)$ and strictly
increases on $(y_0, y_1)$. We have $L(0)\leq L(y_1)$ because $\beta\in
[1/2, 1)$. By continuity, for any $0<y<y_0$ there exists a~unique
$\tilde{y}\in(y_0, y_1)$ that satisfies\vspace*{-2pt}
\[
L(y)=L(\tilde{y}),
\]
which reduces to (\ref{map}) after routine algebra.
\end{pf*}

\begin{pf*}{Proof of Claim \ref{claim3}}
We only prove (\ref{uv1}); the proof of (\ref{uv2}) is similar.
For $0<y<y_0,$ define
\[
D(y) =L\bigl((\beta^{-1}-1)(y_1-y)\bigr)-L(y).
\]
Direct calculation using (\ref{ly}) gives
\begin{eqnarray*}
D(y)&=&(1-\beta)^{p/q}\bigl[\bigl(y_1-(\beta^{-1}-1)(y_1-y)\bigr)^p
\\
&&\hspace*{45pt}{}-(\beta^{-1}-1) \bigl((1-\beta)^{-1} \beta y\bigr)^p -(2-\beta
^{-1})(y_1-y)^p\bigr]
\\
&\leq& 0,
\end{eqnarray*}
where the inequality follows from Jensen's inequality
\[
\bigl(\alpha u+(1-\alpha)v\bigr)^p\leq\alpha u^p+(1-\alpha) v^p,\qquad p>1,
\]
with
\[
\alpha=\beta^{-1}-1,\qquad u=\frac{\beta y}{1-\beta},\qquad v=y_1-y.
\]
That is,
%
\begin{equation}
\label{keyineq}
L\bigl((\beta^{-1}-1)(y_1-y)\bigr)\leq L(y)=L(\tilde{y}),\qquad 0<y<y_0.
\end{equation}
By the strict monotonicity of $L(\cdot)$ on the interval $(y_0, y_1)$,
if $(\beta^{-1}-1)(y_1-y)> \tilde{y}$, then $L((\beta^{-1}-1)(y_1-y))>
L(\tilde{y})$, which contradicts (\ref{keyineq}). Hence
\[
(\beta^{-1}-1)(y_1-y)\leq\tilde{y},
\]
as required.
\end{pf*}

To prove Claim \ref{claim2}, we use Proposition \ref{lem}. Define
\[
Q_\alpha(u, v)=\cases{
\dfrac{u^{\alpha}-v^\alpha}{u-v}, &\quad $u, v>0, u\neq v$,
\cr
\alpha u^{\alpha-1}, &\quad $u=v>0$.
}
\]

\begin{proposition}
\label{lem}
If $0<\alpha\leq1$, then $Q_\alpha(u, v)$ decreases in each of $u,
v>0$; if $\alpha>1$, then $Q_\alpha(u, v)$ increases in each of $u, v>0$.
\end{proposition}

Proposition \ref{lem} follows from basic properties of the generalized
logarithmic mean (\cite{2003Bullen}, pages~386--387).

\begin{pf*}{Proof of Claim \ref{claim2}}
For $0<y<y_0,$ define
\[
u=\beta y,\qquad v=(1-\beta)(y_1-y),\qquad\tilde{u}=\beta\tilde
{y},\qquad\tilde{v}=(1-\beta)(y_1-\tilde{y}).
\]
Claim \ref{claim3} says that $v\leq\tilde{u}$ and $u\leq\tilde{v}$.
Applying Proposition \ref{lem}, we obtain
%
\begin{equation}
\label{m1}
Q_{p/q}(u, v)\geq Q_{p/q}(\tilde{v}, \tilde{u}),\qquad 1<p\leq2,
\end{equation}
and
%
\begin{equation}
\label{m2}
Q_{p/q}(u^{-1}, v^{-1})\geq Q_{p/q}(\tilde{v}^{-1}, \tilde
{u}^{-1}),\qquad p>2.
\end{equation}
After routine algebra, (\ref{m1}) (for $1<p\leq2$) and (\ref{m2}) (for
$p>2$) reduce to (\ref{ineq3}).
\end{pf*}

\end{appendix}

\section*{Acknowledgements}
The author would like to thank Taizhong Hu, Moshe Shaked and Yosef
Rinott for stimulating discussions.

\printhistory

\end{document}